\documentclass[12pt]{amsart}
\usepackage{amsfonts,amssymb,amsthm,eucal}

\newtheorem{thm}{Theorem}
\newtheorem{cor}[thm]{Corollary}
\newtheorem{lemma}[thm]{Lemma}
\newtheorem{prop}[thm]{Proposition}

\newcommand{\R}{\mathbb{R}}

\newcommand{\inprod}[2]{\left\langle #1, #2 \right\rangle}

\newcommand{\vr}{\mathrm{vr}}

\newcommand{\Ent}{\mathrm{Ent}}
\newcommand{\Var}{\mathrm{Var}}

\allowdisplaybreaks

\author[M. Meckes]{Mark W. Meckes}
\email{mark@math.stanford.edu}
\address{Department of Mathematics, Stanford University, Stanford,
California 94305, U.S.A.}

\title[Transportation inequalities]{Some remarks on transportation
 cost and related inequalities}

\begin{document}

\begin{abstract}
We discuss transportation cost inequalities for uniform measures on
convex bodies, and connections with other geometric and functional
inequalities. In particular, we show how transportation inequalities
can be applied to the slicing problem, and prove a new log-Sobolev-type
inequality for bounded domains in $\R^n$.
\end{abstract}

\maketitle

\section{Introduction}\label{S:Intro}

We work in $\R^n$ equipped with its standard inner product
$\inprod{\cdot}{\cdot}$ and Euclidean norm $|\cdot|$. $|A|$ also denotes
the volume (Lebesgue measure) of a measurable set $A$. $D_n$ is 
the Euclidean ball of volume one. For a measurable set
$A$ with $0 < |A| < \infty$, $m_A$ denotes the uniform probability measure
on $A$, that is, $m_A(B) = \frac{|A\cap B|}{|A|}$. The symbols $\mu$ and 
$\nu$ will always stand for Borel probability measures on $\R^n$.

We first introduce two different ways to quantify the difference between
two probability measures. First, for $p\ge 1$, the ($L_p$) {\em Wasserstein
distance} between $\mu$ and $\nu$ is
$$W_p(\mu,\nu) = \inf_\pi \left(\int |x-y|^p \ d\pi(x,y)\right)^{1/p},$$
where $\pi$ runs over probability measures on $\R^n \times \R^n$ with
marginals $\mu$ and $\nu$. We will be interested mainly in the special cases
$p=1,2$. Second, if $\nu \ll \mu$, the {\em relative entropy} of $\nu$
with respect to $\mu$ is
$$H(\nu | \mu) = \int \log \left(\frac{d\nu}{d\mu}\right) \ d\nu.$$
The ($L_p$) {\em transportation cost constant} $\tau_p(\mu)$ is the
largest constant $\tau$ such that
\begin{equation}\label{E:TCI}
W_p(\mu, \nu) \le \sqrt{\frac{2}{\tau}H(\nu|\mu)}
\end{equation}
for every $\nu \ll \mu$. An inequality of the form of (\ref{E:TCI}) is
referred to as a {\em transportation cost inequality} for $\mu$. Note that
if $p \le q$, then $W_p \le W_q$ by H\"older's inequality, and hence 
$\tau_p(\mu) \ge \tau_q(\mu)$.

Transportation cost inequalities are by now well known as a method to derive
measure concentration (cf.\ \cite[Chapter 6]{Ledoux}). 
In fact, as follows from Bobkov and G\"otze's 
dual characterization of the $L_1$ transportation cost inequality \cite{BG},
$\tau_1(\mu)$ is equivalent to the best constant $\alpha$ in the
normal concentration inequality:
\begin{equation}\label{E:conc}
\mu (\{ x \in \R^n : |F(x)| \ge t \})
  \le 2e^{-\alpha t^2} \mbox{ for $t>0$}
\end{equation}
for all $1$-Lipschitz functions $F$ with $\int F d\mu = 0$.

In this paper we consider transportation cost inequalities for uniform
measures on convex bodies. In the next section we show that such inequalities
can be applied to the slicing problem. In the last section we discuss their
relationship with Sobolev-type functional inequalities, and present a
logarithmic Sobolev inequality with trace for bounded domains in $\R^n$.

\section{Relation to the slicing problem}

We recall the following definitions and facts about isotropic convex bodies
(see \cite{MP}).
A convex body $K$ is called {\em isotropic} if
\begin{enumerate}
\item its centroid is $0$,
\item $|K|=1$, and
\item there is a constant $L_K > 0$ such that
  $$\int_K \inprod{x}{y}^2 \ dx = L_K^2 |y|^2$$
  for all $y\in \R^n$.
\end{enumerate}
Every convex body $K$ has an affine image $T(K)$ (unique up to orthogonal
transformations) which is isotropic; the {\em isotropic constant} of $K$
is defined as $L_K = L_{T(K)}$. The isotropic constant also has the extremal
characterization
\begin{equation}\label{E:Ext}
L_K = \min_T \left(\frac{1}{n|K|^{1+2/n}}\int_{T(K)} |x|^2 \ dx\right)^{1/2},
\end{equation}
where $T$ runs over volume-preserving affine transformations of $\R^n$,
with equality iff $K$ is isotropic.
The slicing problem for convex bodies asks whether there is a universal
constant $c$ such that $L_K \le c$ for all convex bodies $K$; see \cite{MP}
for extensive discussion and alternate formulations.

If $K, B \subset \R^n$ are convex bodies, the {\em volume ratio} of $K$ in
$B$ is
$$\vr(B,K) = \min_T \left(\frac{|B|}{|T(K)|}\right)^{1/n},$$
where $T$ runs over affine transformations of $\R^n$ such that $T(K) \subset
B$. The following lemma indicates the relevance of transportation cost
inequalities to the slicing problem.

\begin{lemma}\label{T:Comparison}
Let $K, B \subset \R^n$ be convex bodies, with $B$ isotropic. Then
$$L_K \le c\bigl(1+\sqrt{\log v}\bigr)v\ \tau^{-1/2},$$
where $\tau = \tau_1(B)$, $v=\vr(B,K)$, and $c$ is an absolute constant.
\end{lemma}
\begin{proof}
We may assume that $K\subset B$ and $|K| = v^{-n}$. If $\delta_0$
denotes the point mass at $0\in \R^n$, then by the triangle inequality
for $W_1$,
\begin{eqnarray*}
\frac{1}{|K|}\int_K |x| \ dx &=& W_1(m_K,\delta_0)
    \le W_1(m_K, m_B) + W_1(m_B,\delta_0) \\
&\le& \sqrt{\frac{2}{\tau}H(m_K|m_B)} + \int_B |x| \ dx \\
&\le& \sqrt{\frac{2}{\tau}\log\frac{1}{|K|}}
    + \left(\int_B |x|^2 \ dx\right)^{1/2} \\
&=& \sqrt{\frac{2n}{\tau}\log v} + \sqrt{n}L_B.
\end{eqnarray*}
Now by applying (\ref{E:conc}) to a linear functional, $L_B 
\le c \ \tau^{-1/2}$. 
On the other hand, by Borell's lemma (see e.g.~\cite[Section 2.2]{Ledoux}), 
there is an absolute constant $c$ such that
$$\left(\frac{1}{|K|}\int_K |x|^2 \ dx\right)^{1/2}
    \le c \frac{1}{|K|}\int_K |x| \ dx.$$
The claim now follows from the extremal characterization of $L_K$
(\ref{E:Ext}).
\end{proof}

An analogous estimate with $\tau = \tau_2(B)$ can be proved more
directly, without Borell's lemma.

In light of the equivalence of $L_1$ transportation cost inequalities and
normal concentration, Lemma \ref{T:Comparison} can also be thought of as
an application of measure concentration to the slicing problem.
Since the Euclidean ball is well known to have normal
concentration, as an immediate corollary we 
obtain the following known fact.

\begin{cor}
If $\vr(D_n, K) \le c$, then $L_K \le c'$, where $c'$ depends only on $c$.
\end{cor}

Recently, Klartag \cite{Klartag} introduced the following isomorphic 
version of the slicing problem: given a convex body $K$, is there a convex 
body $B$ such that $L_B \le c_1$ and $d(B,K) \le c_2$, where $d$ is
Banach-Mazur distance? In the case that $K$ and $B$ are 
centrally symmetric, Klartag solved this problem in the affirmative, up 
to a logarithmic (in $n$) factor in $c_2$. 
Lemma \ref{T:Comparison} suggests
approaching the slicing problem via a modified version of the isomorphic
problem: given a convex body $K$, can one find a ``similar'' body $B$
such that $\tau_1(B)$ is large when $B$ is in isotropic
position? Notice that while
Klartag's result uses Banach-Mazur distance to quantify ``similarity'' of
bodies, in Lemma \ref{T:Comparison} it is the weaker measure 
of volume ratio which is relevant. It also seems that this approach via
transportation cost is less sensitive to central symmetry
than more traditional methods of asymptotic convexity.

Finally, we remark that the real point of the proof of Lemma 
\ref{T:Comparison} is that moments of the Euclidean norm on convex bodies,
thought of as functionals of the bodies,
are Lipschitz with respect to Wasserstein distances on uniform measures. 
This suggests an alternative approach to the slicing
problem, related to the
one discussed above, of directly studying optimal (or near-optimal)
probability measures $\pi$ in the definition of $W_p(m_K,m_B)$ for $p=1,2$.
Particularly in the case $p=2$ a great deal is known about the optimal
$\pi$; see \cite{Villani} for an excellent survey.

\section{Functional inequalities}\label{S:Inequalities}

The {\em entropy} of $f:\R^n \to \R_+$ with respect to $\mu$ is
$$\Ent_\mu (f) = \int f \log \left(\frac{f}{\int f d\mu}\right) d\mu,$$
and the {\em variance} of $f:\R^n \to \R$ with respect to $\mu$ is
$$\Var_\mu (f) = \int f^2 d\mu  - \left(\int f d\mu\right)^2.$$
The {\em logarithmic Sobolev constant} $\rho(\mu)$ is the largest constant
$\rho$ such that
\begin{equation}\label{E:LSI}
\Ent_\mu (f^2) \le \frac{2}{\rho} \int |\nabla f|^2 \ d\mu
\end{equation}
for all smooth $f\in L_2(\mu)$. The {\em
spectral gap} $\lambda(\mu)$ is the largest constant $\lambda$ such that
\begin{equation}\label{E:SGI}
\Var_\mu(f) \le \frac{1}{\lambda} \int |\nabla f|^2 \ d\mu
\end{equation}
for all smooth $f\in L_2(\mu)$.
It is well known (cf.\ \cite{Ledoux}) that a logarithmic Sobolev inequality
for $\mu$ implies normal concentration (and hence an $L_1$ transportation
cost inequality, by Bobkov and G\"otze's result \cite{BG}) and a spectral
gap inequality implies exponential concentration. A result of Otto and 
Villani \cite{OV} shows further that
$$\rho(\mu) \le \tau_2(\mu) \le \lambda(\mu)$$
for any absolutely continuous $\mu$. Thus transportation cost inequalities
are somehow intermediate between these Sobolev-type functional inequalities,
and it is of interest here to consider what is known about $\rho(K)$
and $\lambda(K)$ for a convex body $K$. We briefly review known results.

Kannan, Lov\'asz, and Simonovits \cite{KLS} showed that
$$\lambda(K) \ge c \left(\frac{1}{|K|}\int_K |x-z|^2 \ dx\right)^{-1},$$
where $z$ is the centroid of $K$.
It is easy to see that this is an optimal estimate in general. By testing
(\ref{E:SGI}) on linear functionals, one can see that
$\lambda(\mu) \le \alpha_1^{-1}$ for any $\mu$, where $\alpha_1$ is the
largest eigenvalue
of the covariance matrix of $\mu$. If $\alpha_1$ is much larger than the
remaining eigenvalues (i.e., $\mu$ is close to being one-dimensional),
then $\int |x|^2 \ d\mu(x) \approx \alpha_1$.
However, this situation is far from isotropicity (in which all the
eigenvalues are equal), and the authors of \cite{KLS} conjecture that
when $K$ is isotropic,
\begin{equation} \label{E:Conj}
\lambda(K) \ge c n \left(\int_K |x|^2 \ dx\right)^{-1}
    = \frac{c}{L_K^2}.
\end{equation}

Bobkov \cite{Bobkov1} estimated $\rho(K)$ in terms
of the $L_{\psi_2}(m_K)$ norm of $|\cdot|$; in the case that $K$ is 
isotropic, this can be combined with a result of Alesker \cite{Alesker}
to yield
\begin{equation}\label{E:B+A}
\rho(K) \ge \frac{c}{n L_K^2}.
\end{equation}
The estimate for $\tau_1(K)$ which follows from (\ref{E:B+A})
also follows by combining Alesker's result with an $L_1$
transportation cost inequality proved recently by Bolley and Villani
\cite{BV} in an extremely general setting.
The estimate (\ref{E:B+A}) misses the level of (\ref{E:Conj})
by a factor of $n$, but in this case the estimate cannot be
sharpened even when $K$ is isotropic: if $K$ is taken to be
the $\ell_1^n$ unit ball, renormalized to have volume one, then exponential
concentration correctly describes the behavior of a linear functional in a
coordinate direction; it can in fact be shown that $\tau_1(K) \approx 
\frac{1}{n}$ in this case. However, in two concrete cases we have best 
possible estimates:
$$\rho(Q_n) \ge c \quad\mbox{and}\quad \rho(D_n) \ge c,$$
where $Q_n$ is a cube of volume 1. 
The estimate for $Q_n$ is probably folklore;
the estimate for $D_n$ is due to Bobkov and Ledoux \cite{BL}.

Finally, we present the following ``doubly homogeneous $L_p$ trace 
logarithmic Sobolev inequality'' for uniform measures on bounded domains,
inspired both by the search for good estimates on $\rho(K)$ for isotropic
$K$, and by the recent work \cite{MV} by Maggi and Villani on trace
Sobolev inequalities.
This seems not to be directly comparable to the classical logarithmic Sobolev
inequality (\ref{E:LSI}), but interestingly is completely insensitive to 
isotropicity or even convexity of the domain. 
For $p>1$ we denote by $q$ the conjugate
exponent $q = \frac{p}{p-1}$, and $\omega_n = \frac{\pi^{n/2}}{\Gamma
(1+n/2)}$ is the volume of the Euclidean unit ball.

\begin{prop}\label{T:TLSI}
Let $\Omega \subset \R^n$ be open and bounded with locally Lipschitz
boundary and let $p\ge 1$. Then
$$\Ent_\Omega(|f|^p) \le \left(\frac{p-1}{n+q}\right)^{p-1}
    \frac{1}{\omega_n^{p/n}|\Omega|^{1-p/n}} \int_\Omega |\nabla f|^p
    + \frac{1}{\omega_n^{1/n} |\Omega|^{1-1/n}} \int_{\partial \Omega}
    |f|^p$$
for every locally Lipschitz $f:\overline{\Omega} \to \R$, where
$\left(\frac{p-1}{n+q}\right)^{p-1}$ is interpreted as $1$ if
$p=1$, and the integral over $\partial \Omega$ is with respect to 
$(n-1)$-dimensional Hausdorff measure.
\end{prop}

In the case $p=2$ and $f|_{\partial \Omega}=0$, we obtain
$$\Ent_\Omega(f^2) \le \frac{1}{(n+2) \omega_n^{2/n} |\Omega|^{1-2/n}}
    \int_\Omega |\nabla f|^2
    \le \frac{c}{|\Omega|^{1-2/n}} \int_\Omega |\nabla f|^2.$$
Notice that
$$\frac{|\Omega|^{2/n}}{(n+2) \omega_n^{2/n}}\
    = \frac{|\Omega|^{2/n}}{n \omega_n^{1+2/n}} \int_{D_n} |x|^2 \ dx
    \le \frac{1}{n|\Omega|} \int_\Omega |x|^2 \ dx$$
with equality only if $\Omega$ is a Euclidean ball. Therefore if one
restricts the logarithmic Sobolev inequality (\ref{E:LSI}) for $\mu=m_\Omega$
to functions
which vanish on the boundary of $\Omega$, one can improve the constant
$\rho$ to $2 (n+2) \omega_n^{2/n} |\Omega|^{-2/n}$, which is always 
stronger by at least a factor of $2$ than the best possible result for 
general $f$, and much stronger still in many cases.

\begin{proof}[Proof of Proposition \ref{T:TLSI}]
The proof is based on the results of Brenier and McCann on mass
transportation via a convex gradient; we refer to \cite{Villani} for
details and references. To begin, we assume that $p > 1$; the case $p=1$
follows the same lines and is slightly simpler. We also assume that $f$ 
is smooth and nonnegative,
$$\frac{1}{|\Omega|} \int_\Omega f^p = 1,$$
and
$$|\Omega| = \left(\frac{n+q}{p-1}\right)^{n/q} \omega_n.$$
We will use the fact that there is a convex function
$\varphi$ such that $\nabla \varphi$ (the {\em Brenier map})
transports the probability measure
$f^p dm_\Omega$ to the probability measure $m_{B_R}$, where $B_R = 
|\Omega|^{1/n} D_n$
is the Euclidean ball normalized so that $|B_R| = |\Omega|$.

By the results of McCann, the Monge-Amp\`ere equation
$$f^p(x) = \det H_A \varphi(x)$$
is satisfied $f^p dm_\Omega$-a.e., where $H_A \varphi$ is the Aleksandrov
Hessian of $\varphi$ (i.e., the absolutely continuous part of the
distributional Hessian $H\varphi$). Using the fact that $\varphi$ is
convex and $\log t \le t-1$ for $t>0$,
$$\log f^p(x) = \log \det H_A\varphi(x) \le \Delta_A\varphi(x) - n,$$
where $\Delta_A\varphi$ is the Aleksandrov Laplacian of $\varphi$ (i.e.,
the trace of $H_A\varphi$). Integrating with respect to $f^p dm_\Omega$
yields
\begin{equation} \label{E:TLSI1}
\frac{1}{|\Omega|}\int_\Omega f^p \log f^p \le \frac{1}{|\Omega|}
    \int_\Omega f^p \Delta_A \varphi - n
    \le \frac{1}{|\Omega|} \int_\Omega f^p \Delta \varphi - n,
\end{equation}
since $\Delta_A \varphi \le \Delta \varphi$ as distributions, where
$\Delta \varphi$ is the distributional Hessian of $\varphi$. Integrating
by parts (cf.\ \cite{MV} for a detailed justification),
\begin{eqnarray}
\frac{1}{|\Omega|} \int_\Omega f^p \Delta \varphi &=&
    -\frac{1}{|\Omega|}\int_\Omega \inprod{\nabla \varphi}{\nabla (f^p)}
    +\frac{1}{|\Omega|}\int_{\partial \Omega}
    \inprod{\nabla \varphi}{\sigma} f^p \nonumber \\
&=& -\frac{p}{|\Omega|}\int_\Omega f^{p-1}
    \inprod{\nabla \varphi}{\nabla f} + \frac{1}{|\Omega|}\int_{\partial \Omega}
    \inprod{\nabla \varphi}{\sigma} f^p, \label{E:TLSI2}
\end{eqnarray}
where $\sigma$ is the outer unit normal vector to $\partial \Omega$.

Now
\begin{equation} \label{E:TLSI3}
\frac{1}{|\Omega|}\int_{\partial \Omega}
    \inprod{\nabla \varphi}{\sigma} f^p
\le \frac{R}{|\Omega|} \int_{\partial \Omega} f^p
    = \frac{1}{\omega_n^{1/n} |\Omega|^{1-1/n}} \int_{\partial \Omega}
    f^p.
\end{equation}
On the other hand, by H\"older's inequality, the definition of mass
transport, and the arithmetic-geometric means inequality,
\begin{eqnarray}
-\frac{p}{|\Omega|}\int_\Omega f^{p-1} \inprod{\nabla \varphi}{\nabla f}
&\le& p \left(\frac{1}{|\Omega|} \int_\Omega f^p |\nabla \varphi|^q
    \right)^{1/q} \left(\frac{1}{|\Omega|} \int_\Omega |\nabla f|^p
    \right)^{1/p} \nonumber \\
&\le& \frac{p-1}{|B_R|} \int_{B_R} |x|^q dx
    + \frac{1}{|\Omega|} \int_\Omega |\nabla f|^p \nonumber \\
&=& \frac{(p-1)n}{n+q}R^q + \frac{1}{|\Omega|} \int_\Omega |\nabla f|^p
\nonumber \\
&=& n + \left(\frac{p-1}{n+q}\right)^{p-1} \frac{1}{\omega_n^{p/n}
    |\Omega|^{1-p/n}} \int_\Omega |\nabla f|^p. \label{E:TLSI4}
\end{eqnarray}

Combining (\ref{E:TLSI1}), (\ref{E:TLSI2}), (\ref{E:TLSI3}), and 
(\ref{E:TLSI4}) yields
$$\Ent_\Omega(f^p) \le \left(\frac{p-1}{n+q}\right)^{p-1}
    \frac{1}{\omega_n^{p/n} |\Omega|^{1-p/n}} \int_\Omega |\nabla f|^p
    + \frac{1}{\omega_n^{1/n} |\Omega|^{1-1/n}} \int_{\partial \Omega}
    f^p.$$

\medskip

Both sides of this inequality have the same homogeneity with respect to
both $f$ and $\Omega$, so the claim follows for general $f$ and $\Omega$
by rescaling, approximation, and the fact that $|\nabla |f|| = |\nabla f|$
a.e.
\end{proof}

\section*{Acknowledgement}
The author thanks S.\ Bobkov for his comments on an earlier version of this
paper.

\bibliographystyle{plain}
\bibliography{transport}

\end{document}